\documentclass[11pt]{article}

\usepackage[margin=1in]{geometry}
\usepackage{amsmath,amssymb,amsthm,mathrsfs,bm}
\usepackage{booktabs}
\usepackage{graphicx}
\usepackage{enumitem}
\usepackage{hyperref}
\usepackage{setspace}
\usepackage{caption}
\usepackage{subcaption}

\hypersetup{
    colorlinks=true,
    linkcolor=blue,
    citecolor=blue,
    urlcolor=blue
}

\title{\textbf{Housing Decisions under Mobility Risk: \\ 
A Stochastic Threshold Approach}}

\author{
Hui Wu\\
Department of Mathematical Sciences\\
Clark Atlanta University\\
\texttt{hwu@cau.edu}
}

\date{April 2026}

\newtheorem{theorem}{Theorem}
\newtheorem{remark}{Remark}

\begin{document}

\maketitle

\begin{abstract}
Housing tenure decisions are often evaluated using simple metrics such as the price-to-rent ratio. However, in markets where households face significant mobility risk, such as those near military installations, these metrics may provide misleading guidance. This paper develops a tractable framework to analyze the rent-versus-buy decision under stochastic price and rent dynamics combined with an uncertain relocation horizon.

We formulate the problem as an optimal stopping problem with a free-boundary structure, yielding a threshold rule for the price-to-rent ratio. While the underlying stochastic control framework follows standard formulations, we apply it to housing markets characterized by persistent mobility shocks and heterogeneous volatility conditions.

The model shows that mobility risk lowers the optimal buying threshold and increases the value of waiting, implying that households may rationally choose to rent even when conventional indicators suggest buying. Furthermore, identical price-to-rent ratios can imply different optimal decisions across markets once relocation risk and volatility are taken into account.

We illustrate these mechanisms using calibrated examples from representative housing submarkets, including metropolitan and military-adjacent markets. The results highlight that standard valuation rules may be systematically biased in high-mobility environments, underscoring the importance of incorporating relocation risk into housing decision frameworks.
\end{abstract}


\section{Introduction}

Housing tenure decisions are fundamentally shaped by uncertainty about future housing costs and the duration of occupancy. In standard settings, households compare the relative costs and benefits of renting versus owning, often using simple metrics such as the price-to-rent ratio. However, in many real-world environments, the decision is complicated by mobility risk—the possibility that households may be forced or induced to relocate at uncertain times. This feature is especially pronounced in housing markets near military installations, where frequent reassignment, deployment, and contract turnover create persistent and stochastic relocation shocks.

In such markets, renting is not merely a short-term consumption choice but also an option that preserves flexibility under uncertainty. A higher likelihood of relocation shortens the effective ownership horizon, reducing the value of committing to a home purchase. At the same time, local housing markets may exhibit substantial heterogeneity in volatility and liquidity, further affecting the trade-off between renting and owning. As a result, identical price-to-rent ratios can imply very different optimal decisions across markets with different mobility risks and volatility conditions.

This paper develops a tractable framework to quantify how mobility risk and market volatility jointly shape optimal housing tenure decisions. We model the buy-versus-rent problem as an optimal stopping problem with a free boundary, in which the household chooses when to switch from renting to owning under stochastic price and rent dynamics and an exogenous relocation hazard. While the underlying stochastic control structure follows a standard formulation, our contribution lies in applying this framework to housing markets with persistent mobility shocks and providing calibrated insights across representative submarkets.

Our analysis delivers several key insights. First, the presence of relocation risk lowers the optimal buying threshold, making households more likely to rent even when traditional valuation metrics would suggest buying. Second, higher market volatility further increases the option value of waiting, reinforcing the tendency to delay homeownership. Third, the interaction between mobility risk and volatility implies that price-to-rent ratios alone are insufficient statistics for housing decisions in high-mobility environments. These findings highlight that conventional rules of thumb may be systematically biased when applied to markets with substantial relocation uncertainty.

To illustrate these mechanisms, we calibrate the model to representative housing submarkets, including metropolitan areas such as Atlanta and military-adjacent markets such as Columbus (Fort Moore). The calibrated results show that differences in mobility risk and volatility can generate meaningful variation in optimal buying thresholds, even when average price-to-rent ratios are similar. For illustration, we also use a linearized value difference representation to visualize the decision boundary and facilitate interpretation.

The remainder of the paper is organized as follows. Section 2 discusses institutional features of military-adjacent housing markets. Section 3 presents the model setup. Section 4 formulates the household’s optimal stopping problem. Section 5 derives the reduced-form representation of the buy value. Section 6 characterizes the Hamilton--Jacobi--Bellman equation and the free-boundary solution. Section 7 provides comparative statics and numerical illustrations. Section 8 concludes.

\section{Why Military-Adjacent Housing Markets Are Different}

Housing markets near military installations exhibit several structural features that distinguish them from conventional metropolitan or suburban housing markets.

\begin{enumerate}[label=\arabic*.]
    \item \textbf{Stable but transient rental demand.} Rental demand is often supported by rotating military households, contractors, and short-term personnel, generating relatively strong occupancy in many submarkets.

    \item \textbf{Uncertain ownership horizon.} Households may be required or encouraged to relocate on short notice, making the holding period for a purchased home highly uncertain.

    \item \textbf{Demand shocks tied to military activity.} Military-related decisions can create local demand shocks, including troop increases, decreases, or mission changes.

    \item \textbf{Possible mismatch between rental yield and resale depth.} Some military-adjacent markets display attractive rent-to-price ratios but weak resale liquidity, especially in smaller regional cities.
\end{enumerate}

These features imply that both the benefits and the risks of homeownership are path-dependent and uncertain.

\section{Model Setup}

Let $P_t$ denote the house price and $R_t$ denote the rental rate at time $t$. We assume that these variables follow correlated stochastic processes:
\begin{equation}
dP_t = \mu_P P_t\,dt + \sigma_P P_t\,dW_t^{(1)},
\end{equation}
\begin{equation}
dR_t = \mu_R R_t\,dt + \sigma_R R_t\,dW_t^{(2)},
\end{equation}
with correlation structure
\begin{equation}
dW_t^{(1)} dW_t^{(2)} = \rho\,dt.
\end{equation}

To reflect military-specific uncertainty, one may allow an additional demand shock process $M_t$ to shift price and rent dynamics:
\begin{equation}
dP_t = (\mu_P + \beta M_t)P_t\,dt + \sigma_P P_t\,dW_t^{(1)},
\end{equation}
\begin{equation}
dR_t = (\mu_R + \alpha M_t)R_t\,dt + \sigma_R R_t\,dW_t^{(2)}.
\end{equation}

For analytical tractability, we focus on the reduced-form state variable
\begin{equation}
X_t = \frac{P_t}{R_t},
\end{equation}
which represents the price-to-rent ratio. Applying It\^o's lemma yields a geometric diffusion of the form
\begin{equation}
dX_t = \mu_X X_t\,dt + \sigma_X X_t\,dW_t,
\end{equation}
where $W_t$ is a standard Brownian motion derived from the correlated processes $W_t^{(1)}$ and $W_t^{(2)}$, and where $\mu_X$ and $\sigma_X$ are functions of the primitive parameters $\mu_P,\mu_R,\sigma_P,\sigma_R,\rho$.

\section{The Household's Problem}

The household initially rents and chooses a stopping time $\tau$ at which to purchase a home. After purchasing, the household continues owning until a random relocation time $T$. We assume that relocation occurs with hazard rate $\lambda$, so that $T$ may be modeled as exponentially distributed.

Let $K$ denote the transaction cost of entering homeownership, including down payment friction, closing costs, and other sunk expenditures. Let $H$ denote the flow value of housing services under ownership, and let $C$ denote the flow cost of ownership, including mortgage financing cost, taxes, insurance, and maintenance.

The value of entering ownership at state $X$ is denoted by $G(X)$. The value of remaining in the rental state is denoted by $V(X)$. The household solves
\begin{equation}
V(X)=\sup_{\tau}\mathbb{E}\left[
-\int_0^\tau e^{-rt}R_t\,dt + e^{-r\tau}G(X_\tau)
\right].
\end{equation}

\section{Buy Value and Reduced Form}

Under reduced-form assumptions, the value of buying can be expressed as
\begin{equation}
G(X) = -K + \frac{H-C}{r+\lambda} + \mathbb{E}\left[e^{-rT}P_T\right].
\end{equation}

With exponentially distributed relocation time and a stylized reduced-form approximation, one may write
\begin{equation}
\mathbb{E}\left[e^{-rT}P_T\right] \approx \frac{\lambda}{r+\lambda}P.
\end{equation}

This term captures the discounted continuation or resale value upon exit. The key trade-off is that buying yields housing-service benefits and possible appreciation, but requires an irreversible commitment under uncertain duration.

\section{HJB Equation and Free Boundary}

In the continuation region, where the household optimally continues renting, the value function satisfies the Hamilton--Jacobi--Bellman equation
\begin{equation}
rV(X) = R + \mu_X X V'(X) + \frac{1}{2}\sigma_X^2 X^2 V''(X).
\end{equation}

A standard candidate solution takes the form
\begin{equation}
V(X) = AX^\gamma + \frac{R}{r},
\end{equation}
where $\gamma$ is the negative root of the characteristic equation
\begin{equation}
\frac{1}{2}\sigma_X^2\gamma(\gamma-1) + \mu_X\gamma - r = 0.
\end{equation}

There exists a critical threshold $X^*$ separating the continuation and stopping regions. At this free boundary, the standard value-matching and smooth-pasting conditions apply:
\begin{equation}
V(X^*) = G(X^*),
\end{equation}
\begin{equation}
V'(X^*) = G'(X^*).
\end{equation}

\begin{theorem}
Under standard regularity conditions, there exists a threshold $X^*$ such that it is optimal to buy if and only if
\[
X_t \le X^*.
\]
Equivalently, when the price-to-rent ratio is sufficiently low relative to the value of waiting, immediate purchase is optimal; otherwise, continued renting is optimal.
\end{theorem}

\begin{remark}
The threshold $X^*$ depends on the interest rate $r$, the volatility $\sigma_X$, the relocation hazard $\lambda$, the transaction cost $K$, and the effective ownership benefit $H-C$.
\end{remark}

\section{Comparative Statics}

The model implies several economically intuitive comparative-static results.

\subsection*{Volatility}
An increase in $\sigma_X$ raises the option value of waiting. Therefore, greater uncertainty in the price-to-rent ratio makes renting more attractive and lowers the buying threshold:
\[
\frac{\partial X^*}{\partial \sigma_X} < 0.
\]

\subsection*{Mobility Risk}
A higher relocation hazard $\lambda$ shortens the expected ownership horizon and reduces the expected payoff from buying. Thus:
\[
\frac{\partial X^*}{\partial \lambda} < 0.
\]

\subsection*{Interest Rate}
A higher interest rate $r$ increases discounting and typically reduces the relative attractiveness of buying:
\[
\frac{\partial X^*}{\partial r} < 0.
\]

These results reflect the same core principle: buying is an irreversible investment, while renting preserves flexibility.

\section{Calibration and Numerical Illustration}

\subsection{Empirical Proxies for Military-Adjacent Markets}

To connect the theoretical model with real housing environments, we use nearby civilian housing submarkets as proxies for military-adjacent conditions.

We consider four representative markets:
\begin{itemize}
    \item \textbf{Atlanta (Dobbins Air Reserve Base):} diversified metropolitan market with strong liquidity
    \item \textbf{Columbus (Fort Moore):} mid-sized city with higher military dependence
    \item \textbf{Fayetteville (Fort Liberty):} highly base-dependent military town
    \item \textbf{San Diego (Naval Base):} high-price coastal market with strong resale demand
\end{itemize}

These markets span a wide range of mobility risk and volatility environments.

\subsection{Parameterization}

We calibrate the model using representative values:

\begin{table}[h!]
\centering
\caption{Illustrative Parameter Values}
\begin{tabular}{lcc}
\toprule
Parameter & Atlanta (Dobbins) & Columbus (Fort Moore) \\
\midrule
$r$ & 5\% & 5\% \\
$\mu_X$ & 1\% & 0.5\% \\
$\sigma_X$ & 0.15 & 0.25 \\
$\lambda$ (mobility) & 0.10 & 0.20 \\
$K$ & 8\% & 8\% \\
\bottomrule
\end{tabular}
\end{table}

\paragraph{Economic Meaning.}
\begin{itemize}
    \item Higher $\sigma_X$: smaller, less diversified markets imply greater volatility.
    \item Higher $\lambda$: more frequent relocation risk.
    \item Same $K$: transaction costs are assumed comparable across the two cases.
\end{itemize}

\subsection{Comparative Statics}

\begin{figure}[h!]
    \centering
    \includegraphics[width=0.75\textwidth]{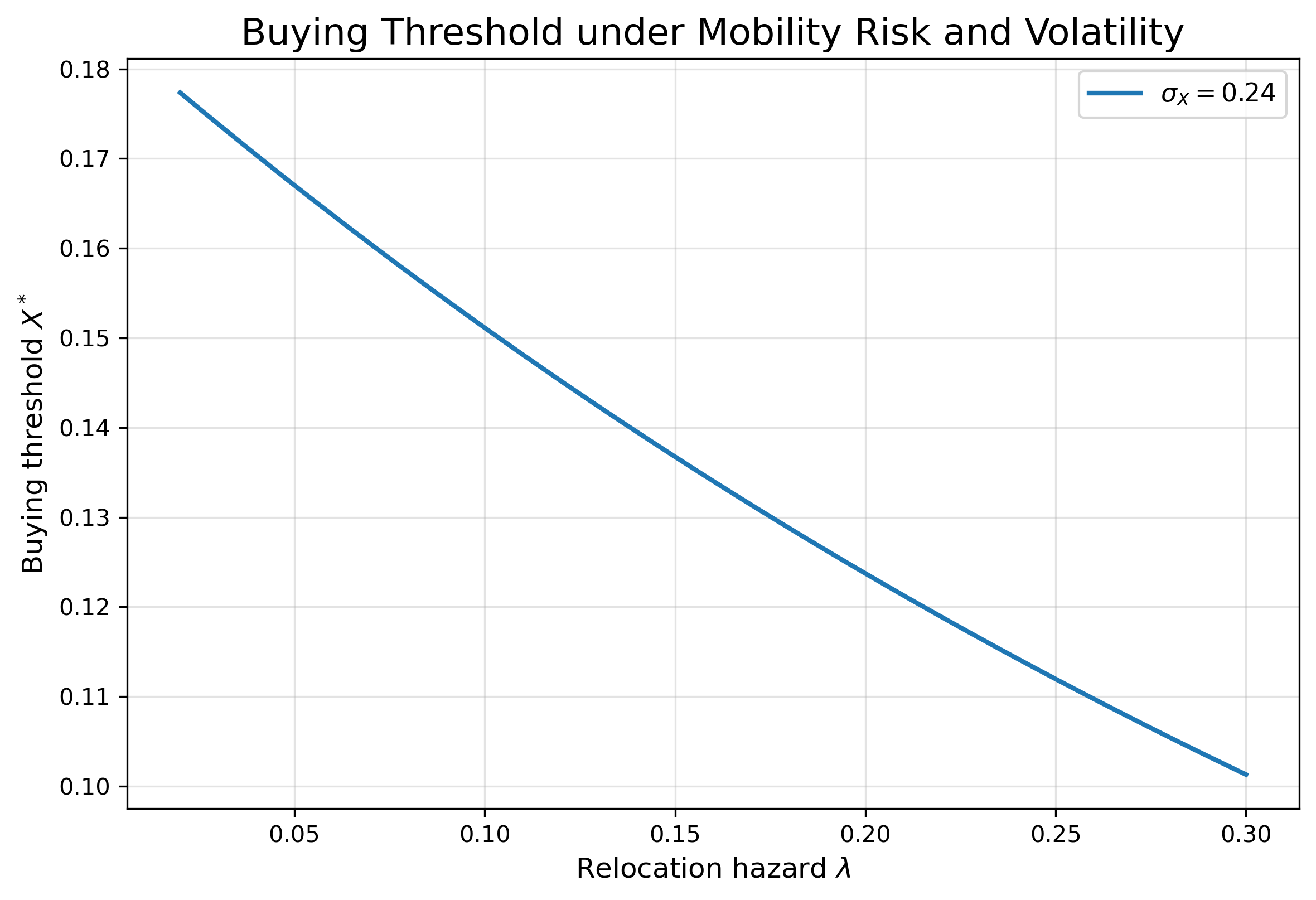}
    \caption{Buying threshold $X^*$ as a function of mobility risk $\lambda$ for different volatility levels.}
    \label{fig:comparative}
\end{figure}

Figure~\ref{fig:comparative} illustrates how the buying threshold $X^*$ varies with mobility risk and volatility.

\paragraph{Interpretation.}
As mobility risk increases, the threshold $X^*$ declines, indicating that households require a lower price-to-rent ratio to justify buying. Higher volatility further lowers the threshold by increasing the option value of waiting.

\subsection{Two-Market Comparison: Atlanta vs Columbus}

\begin{figure}[h!]
    \centering
    \includegraphics[width=0.75\textwidth]{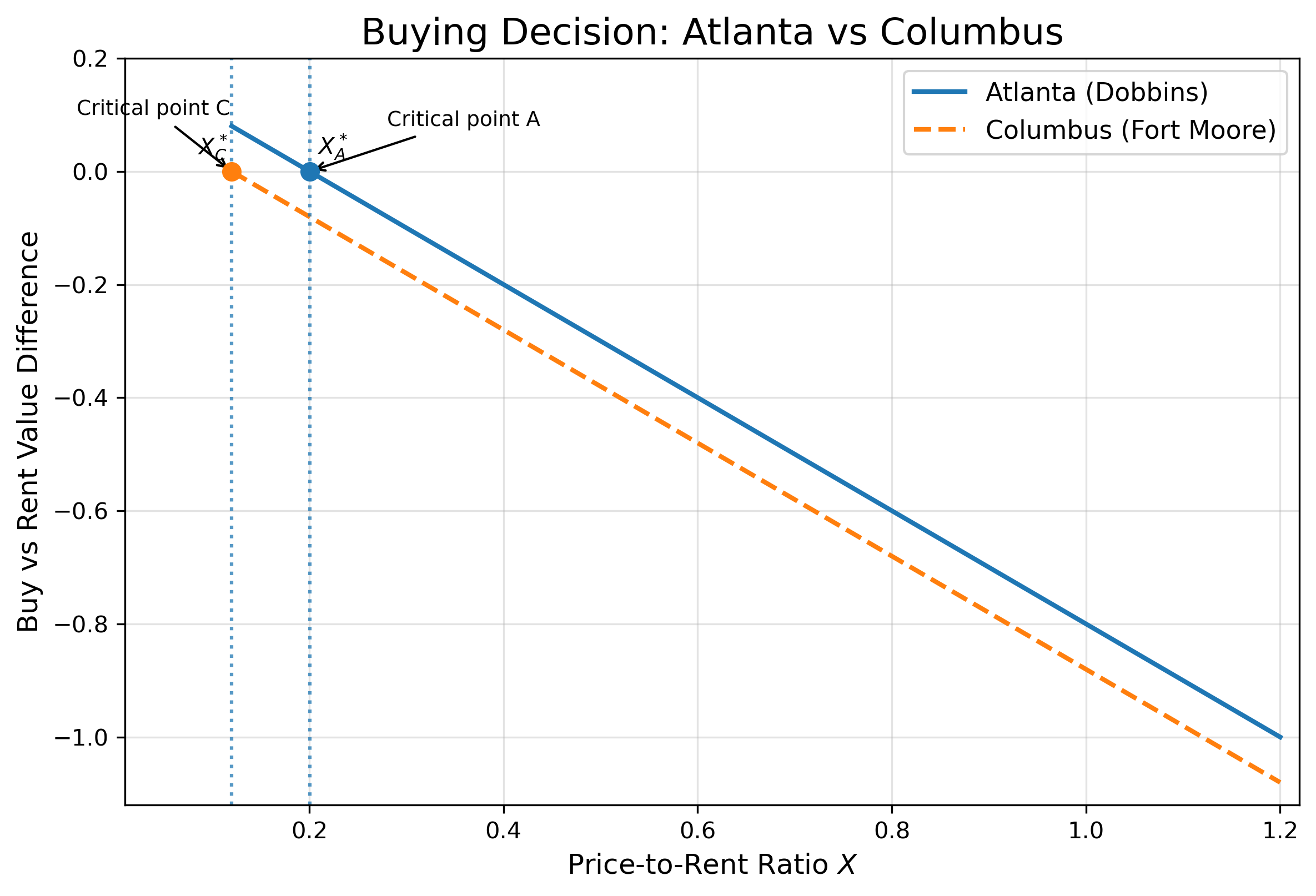}
    \caption{Buying decision comparison between Atlanta and Columbus markets. The vertical dashed lines indicate the threshold values $X^*$ for each market.}
    \label{fig:atl_vs_col}
\end{figure}

Figure~\ref{fig:atl_vs_col} compares the buy-versus-rent decision across two calibrated markets. For illustration, we use a linearized value difference representation to visualize the decision boundary.

\paragraph{Interpretation.}
Atlanta exhibits a higher threshold, indicating that households are more willing to buy. In contrast, Columbus has a lower threshold due to higher mobility risk and volatility, making renting more attractive even when prices appear relatively low. This result highlights that identical price-to-rent ratios can imply different optimal decisions across markets.

\subsection{Threshold Map Across Markets}

\begin{figure}[h!]
    \centering
    \includegraphics[width=0.85\textwidth]{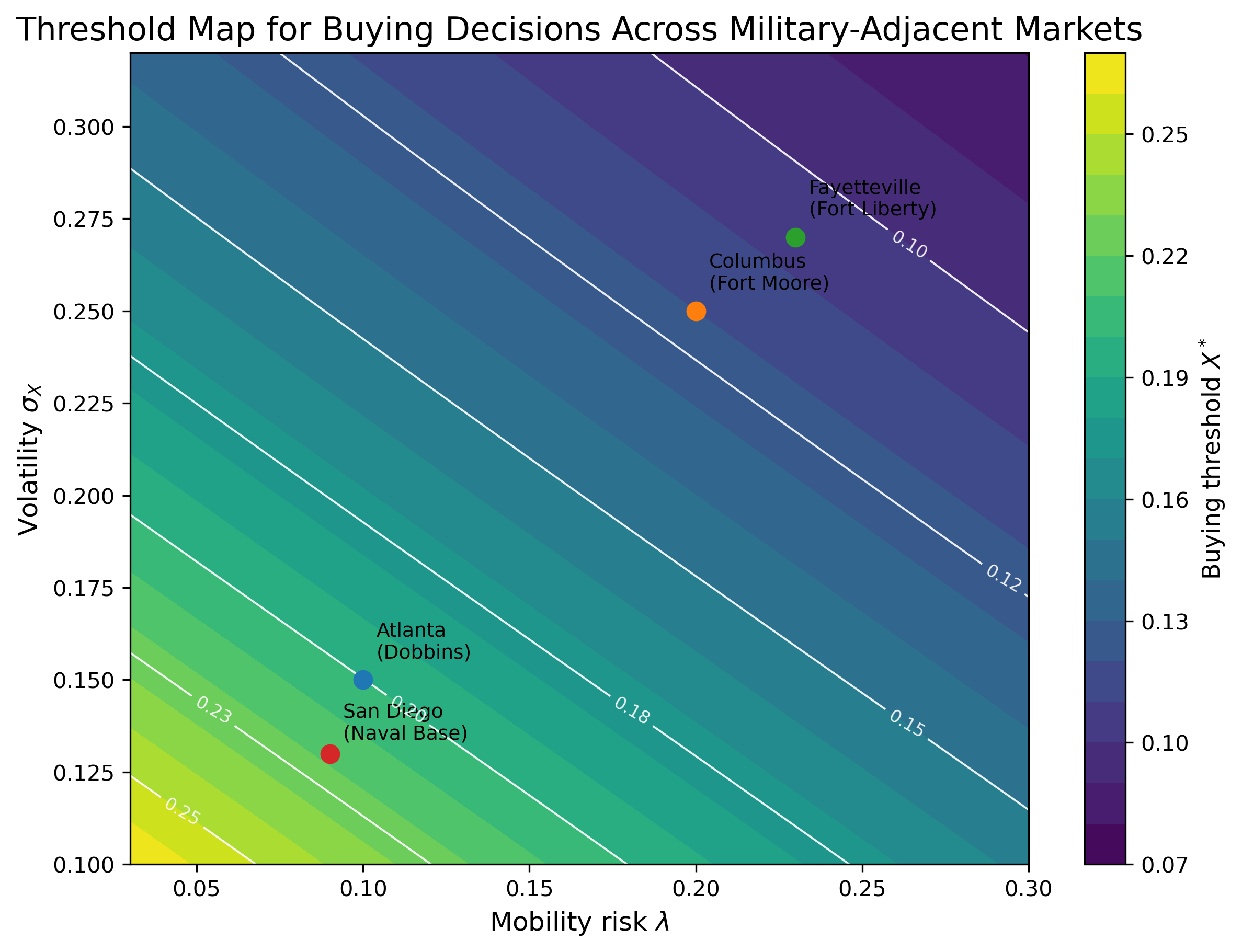}
    \caption{Threshold map for buying decisions across military-adjacent housing markets. The contour levels represent the critical threshold $X^*$ below which buying becomes optimal.}
    \label{fig:map}
\end{figure}

Figure~\ref{fig:map} provides a broader view by mapping the threshold $X^*$ across mobility risk and volatility dimensions.

\paragraph{Interpretation.}
The four cities occupy distinct regions of the parameter space:
\begin{itemize}
    \item Atlanta: moderate mobility, moderate volatility
    \item Columbus: higher mobility and volatility
    \item Fayetteville: high dependence on military demand
    \item San Diego: low volatility despite high prices
\end{itemize}

\paragraph{Key Insight.}
A high rent-to-price ratio does not imply that buying is optimal. In military-adjacent markets, elevated mobility risk and resale uncertainty significantly increase the option value of waiting.

\subsection{Summary}

The numerical results confirm the theoretical predictions:
\begin{itemize}
    \item Buying is an irreversible investment under uncertainty.
    \item Renting provides a real option to delay commitment.
    \item The value of flexibility is particularly large in high-mobility markets.
\end{itemize}

\section{Conclusion}

This paper develops a stochastic framework for the rent-versus-buy decision in housing markets near military installations. Rather than treating the problem as a deterministic comparison of monthly housing costs, we model it as an optimal stopping problem with a free boundary. The resulting threshold policy shows that the decision to buy depends not only on the current price-to-rent ratio, but also on volatility, mobility risk, transaction costs, and the value of flexibility.

The broader economic interpretation is straightforward. Buying is a partially irreversible investment under uncertainty, while renting preserves an option to delay commitment. In military-adjacent markets, where uncertain holding horizons and demand shocks are especially relevant, this option value can be substantial.

Future work may extend the model in several directions, including richer regime-switching demand processes, mortgage-rate dynamics, heterogeneous household preferences, and more detailed empirical estimation using installation-specific panel data.

\section*{References}

\begin{enumerate}
\item Brueckner, J. K. (1997). 
Consumption and investment motives and the portfolio choices of homeowners. 
\textit{Journal of Real Estate Finance and Economics}, 15, 159--180.

\item Capozza, D. R., \& Helsley, R. W. (1990). 
The stochastic city. 
\textit{Journal of Urban Economics}, 28(2), 187--203.

\item Dixit, A. K., \& Pindyck, R. S. (1994). 
\textit{Investment Under Uncertainty}. 
Princeton University Press.

\item Gallagher, J. (2014). 
Learning about an infrequent event: Evidence from flood insurance take-up. 
\textit{American Economic Journal: Applied Economics}, 6(3), 206--233.

\item Genesove, D., \& Mayer, C. (2001). 
Loss aversion and seller behavior: Evidence from the housing market. 
\textit{Quarterly Journal of Economics}, 116(4), 1233--1260.

\item Henderson, J. V., \& Ioannides, Y. M. (1983). 
A model of housing tenure choice. 
\textit{American Economic Review}, 73(1), 98--113.

\item McDonald, R., \& Siegel, D. (1986). 
The value of waiting to invest. 
\textit{Quarterly Journal of Economics}, 101(4), 707--727.

\item Piazzesi, M., Schneider, M., \& Tuzel, S. (2007). 
Housing, consumption, and asset pricing. 
\textit{Journal of Financial Economics}, 83(3), 531--569.

\item Pindyck, R. S. (1991). 
Irreversibility, uncertainty, and investment. 
\textit{Journal of Economic Literature}, 29(3), 1110--1148.

\item Sinai, T., \& Souleles, N. S. (2005). 
Owner-occupied housing as a hedge against rent risk. 
\textit{Quarterly Journal of Economics}, 120(2), 763--789.

\item Titman, S. (1985). 
Urban land prices under uncertainty. 
\textit{American Economic Review}, 75(3), 505--514.
\end{enumerate}

\end{document}